\documentclass[12pt,finalcls,onecolumn]{IEEEtran}
\usepackage{graphicx}
\usepackage{amsmath}
\usepackage{mathrsfs}
\usepackage{float}
\usepackage{verbatim}
\usepackage{enumerate}
\usepackage{layout}

\usepackage{amsmath, amsthm, amsfonts}
\usepackage{amssymb}

\newtheorem{thm}{Theorem}[section]

\newtheorem{defn}[thm]{Definition}
\newtheorem{rem}[thm]{Remark}
\newtheorem{example}[thm]{Example}
\def\RR{\mathbb{R}}

\def\CC{\mathbb{C}}
\def\ualpha{{\underline\alpha}}
\def\oalpha{{\overline\alpha}}
\def\ubeta{{\underline\beta}}
\def\obeta{{\overline\beta}}
\def\oa{{\overline a}}
\def\ua{{\underline a}}
\def\ob{{\overline b}}
\def\ub{{\underline b}}
\def\ba{{\bf a}}

\def\ubb{{\underline\bb}}
\def\obb{{\overline\bb}}

\def\bb{{\bf b}}

\def\Thetasub1{{\Theta_1}}

\def\rb1{{]_1}}

\def\bc{{\bf c}}

\def\bd{{\bf d}}

\def\uy{{\underline  y}}
\def\oy{{\overline  y}}

\def\b0{{\boldsymbol {0}}}

\title{  Algorithm to   Compute a Kharitonov-Type Sector Containing All  Roots of Hurwitz  Interval Polynomials }
\author{David~Hertz%~\IEEEmembership{Senior Member,~IEEE}
\thanks{David~Hertz is an independent researcher, Akko, Israel,\\
e-mail: drdavidhertz@gmail.com. Manuscript version of \today.}}

\begin{document}
\maketitle

\begin{abstract}
%\boldmath

Given a complex interval Hurwitz polynomial, we  present a Kharitonov-type algorithm that uses up to 16 additional Kharitonov polynomials
to check whether or not
the angles  of all  the roots of  a given interval polynomial are in a given sector
$[0.5\pi+\alpha,1.5\pi-\beta], \alpha,\beta\in(0,0.5\pi)$. If  the angles of all the roots of a given  interval polynomial are in a given sector we call it  a {\bf containing sector}.\\
If the given interval polynomial is real  up to 8  additional Kharitonov polynomials are needed to check whether or not
 $[0.5\pi+\alpha,1.5\pi-\alpha], \alpha\in(0,0.5\pi)$ is a containing sector.
Next, by using the bisection algorithm we obtain   in the complex case  the containing sector ${\cal K}:=[0.5\pi+\alpha_K,1.5\pi-\beta_K]$
and in the real case  containing the sector ${\cal K}:=[0.5\pi+\alpha_K,1.5\pi-\alpha_K]$. Notice that the bisection algorithm produces in the real case the interval
$[\alpha_K]:=[\ualpha_K,\oalpha_K]$ and in the complex  case the intervals  $[\alpha_K]:=[\ualpha_K,\oalpha_K]$  and $[\beta_K]:=[\ubeta_K,\obeta_K]$. Hence,    in the real  case ${\cal K}=[0.5\pi+\ualpha_K,1.5\pi-\ualpha_K]$ is a containing sector and in the complex case 
 ${\cal K}=[0.5\pi+\ualpha_K,1.5\pi-\ubeta_K]$ is a containing sector.  
By letting the imaginary parts of the coefficients in the eight Kharitonov polynomials be zero we obtain that the  first four Kharitonov polynomials are  the same as   the remaining four Kharitonov polynomials. Thus, we  could  trivially prove Kharitonov real case from Kharitonov complex case.\\
Simulation results led us to conjecture that the minimal containing sector equals to the containing sector of all the vertex polynomials or  to a  given subset thereof.
 
\end{abstract}

% Note that keywords are not normally used for peer review papers.

% For peer-reviewed papers, you can put extra information on the cover
% \fi
%
% For peer review papers, this IEEEtran command inserts a page break and
\IEEEpeerreviewmaketitle

\section{Introduction}

Kharitonov theorems provide a fundamental framework for
the robust stability analysis of interval polynomials and have motivated
a large body of research on parametric uncertainty and Hurwitz stability.
In this paper, we develop an algorithm to compute a Kharitonov-type
sector containing all roots of a real or complex Hurwitz interval polynomial.

We define a complex interval  polynomial by
\begin{align}
P(s;[\ba+i\bb]):=\sum_{n=0}^N[a_n+ib_n]s^n,
\end{align}
where $[a_n+ib_n]:=[a_n]+i[b_n]=$
$[\ua_n,\oa_n]+i[\ub_n,\ob_n], [\ua_n,\oa_n],[\ub_n,\ob_n]\subset\RR~
(n=0,1,...,N$),\\
 $[\ba]:=([a_0],[a_1],...,[a_N])$, and\\ $[\bb]:=([b_0],[b_1],...,[b_N])$.

We define a real interval  polynomial by
\begin{align}
P(s;[\ba]):=\sum_{n=0}^N[a_n]s^n,
\end{align}
where $ [a_n]:=[\ua_n,\oa_n]\subset\RR ~(n=0,1,...,N$) 
and $[\ba]:=([a_0],[a_1],...,[a_N])$.
Kharitonov  proved   \cite{KharitonovR} that if 4 real vertex polynomials are Hurwitz then  the real interval polynomial 
 $P(s;[\ba])$ is Hurwitz; and,  later on proved   \cite{Kharitonov8} that if 8 complex vertex polynomials are Hurwitz then  the complex interval polynomial 
 $P(s;[\ba+i\bb])$ is Hurwitz. 
 For what follows we abbreviate Kharitonov 4 and 8 polynomials theorem by K4PT and K8PT, respectively.
 The book \cite{RobCont} devotes Chapter~5 to Kharitonov theorems including  proofs by extending the Hermite-Bieler Interlacing Theorem. 
There are other  proofs of K4PT and K8PT, see  the references in \cite{RobCont},\cite{HK}, and {\cite{Bistriz} to mention but a few. The proof by Bistriz {\cite{Bistriz} based   on the evaluation of complex lossless positive-real functions and  
their relation to Hurwitz polynomials may have potential for further  extensions.

The purpose of this article, a corrected and slightly extended version of the recent
preprint \cite{HertzPreprint}, is to present an algorithm that determines
whether both $[0.5\pi+\alpha,1.5\pi)$ and $(0.5\pi,1.5\pi-\beta]$ are
containing sectors for $P(s;[\ba+i\bb])$, where
$\alpha,\beta\in(0,0.5\pi)$. In that case,
$[0.5\pi+\alpha,1.5\pi-\beta]$ is a containing sector for $P(s;[\ba+i\bb])$.

 Next by using the bisection algorithm we arrive at a sub-optimal Kharitonov-type containing sectors
${\cal K}=[0.5\pi+\alpha_K,1.5\pi-\beta_K]$  for  $P(s;[\ba+i\bb])$ and 
${\cal K}=[0.5\pi+\alpha_K,1.5\pi-\alpha_K]$ for   $P(s;[\ba])$.

In Section~\ref{sec:KC} we  present the proposed algorithm to compute ${\cal K}=[0.5\pi+\alpha_K,1.5\pi-\beta_K]$  for\\
  $P(s;[\ba+i\bb])$
and a comprehensive example that led us to a conjecture.
Next, in Section~\ref{sec:KR} we  present the proposed algorithm to compute $[0.5\pi+\alpha_K,1.5\pi-\alpha_K]$  for  $P(s;[\ba])$ and a comprehensive example  that led us to a similar conjecture.  

By letting the imaginary part of the coefficients of the eight Kharitonov polynomials in K8PT to be zero we obtain that the  first four Kharitonov polynomials are the same as  the remaining four Kharitonov polynomials. 
This immediately implies K4PT as a special case of K8PT.

Finally, in  Section~\ref{sec:Conc} we  give the conclusion and areas for future research.

\section{Smallest  Kharitonov Sector - Complex Case}
\label{sec:KC}
Let
\begin{align}
P(s;[\ba+i\bb])=\sum_{n=0}^N[a_n+ib_n]s^n, 0\not \in [a_N+ib_N].
\end{align}
The corresponding eight Kharitonov polynomials are
\begin{align}
\label{eq:K8}
K_1(s)&=(\oa_0+i\ob_0)+(\oa_1+i\ub_1)s+(\ua_2+i\ub_2)s^2+\\ \nonumber
&\quad\,(\ua_3+i\ob_3)s^3+(\oa_4+i\ob_4)s^4+\cdots \\ \nonumber
K_2(s)&=(\oa_0+i\ub_0)+(\ua_1+i\ub_1)s+(\ua_2+i\ob_2)s^2+\\ \nonumber
 &\quad\,(\oa_3+i\ob_3)s^3+(\oa_4+i\ub_4)s^4+\cdots \\ \nonumber
K_3(s)&=(\ua_0+i\ob_0)+(\oa_1+i\ob_1)s+(\oa_2+i\ub_2)s^2+\\ \nonumber
&\quad\,(\ua_3+i\ub_3)s^3+(\ua_4+i\ob_4)s^4+\cdots \\ \nonumber
K_4(s)&=(\ua_0+i\ub_0)+(\ua_1+i\ob_1)s+(\oa_2+i\ob_2)s^2+\\ \nonumber
&\quad\,(\oa_3+i\ub_3)s^3+(\ua_4+i\ub_4)s^4+\cdots \\ \nonumber
K_5(s)&=(\oa_0+i\ob_0)+(\ua_1+i\ob_1)s+(\ua_2+i\ub_2)s^2+\\ \nonumber
&\quad\,(\oa_3+i\ub_3)s^3+(\oa_4+i\ob_4)s^4+\cdots \\ \nonumber
K_6(s)&=(\oa_0+i\ub_0)+(\oa_1+i\ob_1)s+(\ua_2+i\ob_2)s^2+\\ \nonumber
&\quad\,(\ua_3+i\ub_3)s^3+(\oa_4+i\ub_4)s^4+\cdots \\ \nonumber
K_7(s)&=(\ua_0+i\ob_0)+(\ua_1+i\ub_1)s+(\oa_2+i\ub_2)s^2+\\ \nonumber
&\quad\,(\oa_3+i\ob_3)s^3+(\ua_4+i\ob_4)s^4+\cdots \\ \nonumber
K_8(s)&=(\ua_0+i\ub_0)+(\oa_1+i\ub_1)s+(\oa_2+i\ob_2)s^2+\\ \nonumber
&\quad\,(\ua_3+i\ob_3)s^3+(\ua_4+i\ub_4)s^4+\cdots.  \nonumber
\end{align}
We have that 
if $K_1(s),\cdots,K_8(s)$ are Hurwitz then $P(s;[\ba+i\bb])$ is Hurwitz. 
In what follows we assume that $P(s;[\ba+i\bb])$ is Hurwitz, i.e., $P(s;[\ba+i\bb])\ne 0, \Re s\ge 0, \ba+i\bb\in [\ba+i\bb]$. 
\begin{rem}
We have:\\
(i) If $0\in [a_0+ib_0]$ then   0 is a root of $P(s;[\ba+i\bb])$ and  $P(s;[\ba+i\bb])$  is not Hurwitz.\\
(ii) If $0\in [a_N+ib_N]$ then when $a_N+ib_N=0$ we obtain  a degree reduction.  We will  now give an example to
 gain more understanding of   degree reduction. Let $\tilde P(s)=2s+s^2$ and  $P(s)=s^2\tilde P(s^{-1})=1+2s$ denote the polynomial 
whose roots are the reciprocals of the roots of $\tilde P(s)$ and vice versa. Hence, since the roots of $\tilde P(s)$ are $s_1=-2, s_2=0$  the roots of $P(s)$ are $s_1=-1/2,s_2=\infty$. So, when $a_N+ib_N=0$ $P(s;[\ba+i\bb])$ has a root  at infinity, say $s_N=\infty$, and
 its degree is $N-1$. 
 Since the angle of  a  root at infinity can be in the right half plane then for $\epsilon>0$ small enough satisfying  $0<|a_N+ib_N|<\epsilon$ we may obtain
a very large root in the RHP of $s$. E.g., $P(s)=10^{-6}(1+i)s^2+s+1\approx s+1$ gives that $P(s)$ one root in the open  LHP and the other large root in the open RHP.
 Therefore, we require that $0\notin [a_N+ib_N]$.

\end{rem}
We define a sector in the complex $s$-plane by
\begin{align}
[0.5\pi+\alpha,1.5\pi-\beta]:=
\{{\rm angle}(s): s\in\CC,\ \\ \nonumber
 0.5\pi +\alpha\le{\rm angle}(s)\le 1.5\pi -\beta\}   ,\alpha,\beta\in(0,0.5\pi).\nonumber
\end{align}
\begin{defn}
We call the sector $[0.5\pi+\alpha,1.5\pi-\beta]$  a {\bf containing sector} for $P(s;[\ba+i\bb])$ if it contains  all the angles of   its roots. 
\end{defn}
Let  ${\cal O}:=[0.5\pi+\alpha_0,1.5\pi-\beta_0],\alpha_0,\beta_0\in(0,0.5\pi)$ denote the smallest containing sector for $P(s;[\ba+i\bb])$.
Notice that ${\cal O}$ must be closed; otherwise, we could improve it.
In what follows we  present an algorithm based on  Kharitonov 8 Polynomial Theorem (K8PT)  and the bisection algorithm 
to produce a sub-optimal containing sector, say
${\cal K}:=[0.5\pi+\alpha_K,1.5\pi-\beta_K]\supseteq {\cal O}$  for $P(s;[\ba+i\bb])$.\\
First, we  present an algorithm based on   applying K8PT to  two newly derived interval polynomials and thus   arrive at a decision whether 
or not both   $([0.5\pi+\alpha,1.5\pi)$  and
$(0.5\pi,1.5\pi-\beta]$    are containing sectors for $P(s;[\ba+i\bb])$. If this is the case   we conclude that 
 $[0.5\pi+\alpha,1.5\pi-\beta]$ is a
 containing sector  for   $P(s;[\ba+i\bb])$.  Notice that to decide whether or not  $[0.5\pi+\alpha,1.5\pi-\beta]$ is a containing sector we need
 to test that up to  $16=2\times 8$  Kharitonov polynomials are Hurwitz\footnote{If any Kharitonov polynomial is not Hurwitz  then $[0.5\pi+\alpha,1.5\pi-\beta]$ is not a containing sector  and we are done.}.\\
 Next, by  combining  this algorithm and the bisection algorithm we  arrive at the desired suboptimal containing
 sector  ${\cal K} =[0.5\pi+\alpha_K,1.5\pi-\beta_K]$  for $P(s;[\ba+i\bb])$.
\\ 
\\
{\bf Case (i): Algorithm to test whether or not     $[0.5\pi+\alpha,1.5\pi)$   is  a containing sector for $P(s;[\ba+i\bb]), \alpha\in(0,0.5\pi).$}\\

Here, by using the proposed algorithm  we decide whether or not  $[0.5\pi+\alpha,1.5\pi)$  is a containing sector for $P(s;[\ba+i\bb])$.
The  idea is to  rotate clock wise (CW)  all the roots of $P(s;[\ba+i\bb])$  by $\alpha$.  
We thus obtain the polynomial $P_\alpha(s;[\bc+i\bd]):=P(e^{i\alpha}s;[\ba+i\bb])$, where 
 $[\bc+i\bd]:=([a_0+ib_0],[a_1+ib_1] e^{i\alpha},...,[a_N+ib_N] e^{iN\alpha})$. This is so since  if  $Q(s)$ is a polynomial of degree $N$,
 $Q(s_0)=0$, and $S(s):=Q(e^{i\alpha}s)$ then  $S(e^{-i\alpha}s_0)=Q(e^{i\alpha}(e^{i\alpha}s_0))=0$, i.e., $e^{-i\alpha}s_0$
 is a root of $S(s)$.
 Notice that since $P(s;[\ba+i\bb])$ is Hurwitz  no roots  will enter into  $P_\alpha(s;[\bc+i\bd])$ from the closed fourth $s$-quadrant  to the  open third $s$-quadrant during the rotation. 
 %Furthermore, if $P_\alpha(s;[\bc+i\bd])$ is Hurwitz then $P(s;[\ba+i\bb])$ has no roots in $(0.5\pi,0.5\pi+\alpha]$.

The interval coefficients of the rotated polynomial  $P_\alpha(s;[\bc+i\bd])$ for $n=0,1,...,N$  are given by

\begin{align}
[c_n+id_n]&=e^{in\alpha}[a_n+ib_n]\\                                                                   \nonumber
&=(\cos(n\alpha)+i\sin(n\alpha))[a_n+ib_n]\\                                                          \nonumber
&=\cos(n\alpha)[a_n]-\sin(n\alpha)[b_n]+\\ \nonumber
&\quad\; i(\sin(n\alpha)[a_n]+\cos(n\alpha))[b_n]),    \nonumber
\end{align}
where 
\begin{equation}
[c_n]=\cos(n\alpha)[a_n]-\sin(n\alpha)[b_n],
\end{equation}
\begin{equation}
[d_n]=\sin(n\alpha)[a_n]+\cos(n\alpha)[b_n].
\end{equation}
For $[y]=[\uy,\oy]\subset \RR$ we have
\begin{equation}
\label{eq:xy}
x[y]=\begin{cases}
[x\uy,x\oy],{\rm~if~}x\ge0\\
[x\oy,x\uy],{\rm~if~}x<0.
\end{cases}
\end{equation} 
If by applying K8PT  to  $P_\alpha(s;[\bc+i\bd])$ we obtain that  it
is Hurwitz then there were no roots of   $P(s;[\ba+i\bb])$ in the interval $(0.5\pi,0.5\pi+\alpha]$, therefore,  $[0.5\pi+\alpha,1.5\pi)$   is a containing sector   for $P(s;[\ba+i\bb])$. Notice that  $(0.5\pi+\alpha,1.5\pi)$ is a containing sector for  $P(s;[\ba+i\bb])$ but so is also $(0.5\pi+\alpha,1.5\pi]$ . 
\\ 
\\
{\bf Case (ii): Algorithm to test whether or not     $(0.5\pi,1.5\pi-\beta]$   is  a containing sector  for $P(s;[\ba+i\bb]), \beta\in(0,0.5\pi).$}\\
This  case is similar to case (i). Here, we  rotate
  counter clock wise (CCW)  all the roots of $P(s;[\ba+i\bb])$  by $\beta$.  
We thus obtain the polynomial $P_{-\beta}(s;[\bc+i\bd]):=P(e^{-i\beta}s;[\ba+i\bb])$, where 
 $[\bc+i\bd]:=([a_0+ib_0],[a_1+ib_1] e^{-i\beta},...,$ $[a_N+ib_N] e^{-iN\beta})$.\\
If by applying K8PT  to  $P_{-\beta}(s;[\bc+i\bd])$ we obtain that it is  Hurwitz then  there were no roots of   $P(s;[\ba+i\bb])$ in the sector $[1.5\pi-\beta,1.5\pi)$, therefore,  $(0.5\pi,1.5\pi-\beta]$   is also a containing sector   for $P(s;[\ba+i\bb])$. 
\begin{example}
Let
\begin{align*}
P(s;[\ba+i\bb])=&([2.9475,3.1475]+i[- 0.43,- 0.23])+\\
&\quad\;([6.455,6.655]+i[- 0.425,- 0.225])s+\\ \nonumber
       & ([4.4,4.6]+i[-0.15,0.05])s^2+\\
       &\quad\;([0.9,1.1]+i[-0.1,0.1])s^3. \nonumber
\end{align*}
We derived this interval polynomial  by computing the coefficients of the polynomial whose roots are $s_1=-1+0.1i, s_2=-1.5+0.2i, s_3=-2-0.25i$ and then
replaced each coefficient, say $a_n+ib_n$, by $[a_n+ib_n-0.1-0.1i, a_n+ib_n+0.1+0.1i]$. 
The corresponding eight Kharitonov polynomials are 
\begin{align*}
K_1(s)&=(3.1475 - 0.23i)+(6.655-0.425i)s +\\ \nonumber
&\quad\;(4.4 - 0.15i)s^2+(0.9+0.1i)s^3\\ \nonumber
K_2(s)&=(3.1475 - 0.43i)+(6.455-0.425i)s +\\
&\quad\;(4.4 +0.05i)s^2+(1.1+0.1i)s^3\\ \nonumber
K_3(s)&=(2.9475 - 0.23i)+(6.655-0.225i)s +\\
&\quad\;(4.6 - 0.15i)s^2+(0.9-0.1i)s^3\\ \nonumber
K_4(s)&=(2.9475 - 0.43i)+(6.455 - 0.225i)s +\\
&\quad\;(4.6 + 0.05i)s^2+(1.1 - 0.1i)s^3 \\ \nonumber
K_5(s)&=(3.1475 - 0.23i)+(6.455 - 0.225i)s +\\
&\quad\;(4.4 - 0.15i)s^2+(1.1 - 0.1i)s^3\\ \nonumber
K_6(s)&=(3.1475 - 0.43)+(6.655 - 0.225i)s +\\
&\quad\;(4.4 + 0.05i)s^2+(0.9 - 0.1i)s^3\\ \nonumber
K_7(s)&=(2.9475 - 0.23i)+(6.455 - 0.425i)s\\
&\quad\; +(4.6 - 0.15i)s^2+(1.1+0.1i)s^3\\ \nonumber
K_8(s)&=(2.9475 - 0.43i)+(6.655 - 0.425i)s +\\
&\quad\;(4.6 + 0.05i )s^2+(0.9 + 0.1i)s^3.\nonumber
\end{align*}
%$(2.9475 - 0.43i)+(6.455 - 0.425i)s+(4.4 -0.15i )s^2+(0.9 - 0.1i)s^3

These eight Kharitonov polynomials are Hurwitz, therefore, $P(s;[\ba+i\bb])$ is Hurwitz.\\
By combining   the algorithm in Case (i) and  the bisection algorithm we obtain that the containing sector $[0.5\pi+\alpha_K,1.5\pi)$ 
for  $P_{\alpha_K}(s;[\bc+i\bd])$ satisfies
$\alpha_K\in[\ualpha_K ,\oalpha_K]=[0.2527,0.2528]\pi=[45.486^\circ , 45.504^\circ ]$.
Notice that $[0.5\pi+\ualpha_K,1.5\pi)$ is a containing sector for
 $P_{\ualpha_K}(s;[\bc+i\bd])$, whereas  
 $P_{\oalpha_K}(s;[\bc+i\bd])$ is not Hurwitz. \\
By combining   the algorithm in Case (ii) and  the bisection algorithm  we obtain that the containing sector $(0.5\pi,1.5\pi-\beta_K]$ 
for  $P_{-\beta_K}(s;[\bc+i\bd])$ satisfies
$\beta_K\in[\ubeta_K ,\obeta_K]=[0.2725,0.2726]\pi=[49.05^\circ , 49.068^\circ ]$.\\
Notice that $(0.5\pi,1.5\pi-\ubeta_K]$ is a containing sector for  $P_{-\ubeta_K}(s;[\bc+i\bd])$, whereas    $P_{-\obeta_K}(s;[\bc+i\bd])$
is not Hurwitz.\\
Hence, ${\cal K}:=[0.5\pi+\ualpha_K,1.5\pi-\ubeta_K]=[135.486^\circ,222.95^\circ]$ is a containing sector for $P(s;[\ba+i\bb])$.
To check this result we generated a million uniformly distributed polynomials from $[\ba+i\bb]$. These polynomials appeared in the sector
${\cal S}=[141.9432^\circ,  214.1573^\circ]$, where  ${\cal K}\supseteq{\cal S} $. 
We also checked all the 256 vertex polynomials of $P(s;[\ba+i\bb])$ and obtained that they are in the sector
${\cal V}=[140.3779^\circ,  215.679^\circ]$, where ${\cal V}(1)$ corresponds to
 $V_1(s)=(3.1475 - 0.23i)+(6.455 - 0.425i)s+(4.4 +0.05i )s^2+(1.1 - 0.1i)s^3$ and
${\cal V}(2)$ corresponds to   $V_2(s)=(2.9475 - 0.43i)+(6.655 - 0.225i)s+(4.4 -0.15i )s^2+(1.1 - 0.1i)s^3$. 
%Notice that 
% both $V_1(s),V_2(s)$ are not  Kharitonov polynomials of $P(s;[\ba+i\bb])$. 
So, we obtained     that  ${\cal K}\supseteq  {\cal O}
 \supseteq{\cal V}\supseteq {\cal S}$ and interestingly,we observed that ${\cal V}\supseteq
  {\cal S}$.\\
%(2.9475 - 0.43i)+(6.455 - 0.425i)s+(4.4 -0.15i )s^2+(0.9 - 0.1i)s^3
\end{example}
\begin{rem}
Notice that the coefficients of $P_\theta(s,[\bc+i\bd])$ depend on $\theta\in\{\alpha,-\beta\}$.
The reason why we claimed that Kharitonov containing sector $[0.5\pi+\alpha_K,1.5\pi-\beta_K]$ is sub-optimal, 
is because the interval coefficients of $P_\theta(s,[\bc+i\bd])$  are dependent.
\end{rem}
\begin{rem}
\label{rem:congC}
Based on the above example 
we conjecture that the minimal containing sector is attained at  the vertex polynomials of $P(s;[\ba+i\bb])$.
If this is the case we have  to test $2^{2(N+1)}$  vertex polynomials unless we can prove that a smaller predetermined subset thereof suffices.
\end{rem}

\section{Smallest Kharitonov Containing Sector - Real Case}
\label{sec:KR}

Let
\begin{align}
P(s;[\ba]):=\sum_{n=0}^N[a_n]s^n
\end{align}
be   a given real interval polynomial.
Kharitonov \cite{KharitonovR} (see also \cite{RobCont} for a proof) proved that the interval polynomial $P(s;[\ba])$ is Hurwitz if and only if
the  following 4  Kharitonov vertex polynomials are Hurwitz.\\
\begin{align*}
K_1(s)&=\ua_0+\ua_1s+\oa_2s^2+\oa_3s^3+\ua_4s^4+\\
&\quad\ua_4s^5+\oa_6s^6+\oa_7s^7+\cdots, \\ %\nonumber
K_2(s)&=\ua_0+\oa_1s+\oa_2s^2+\ua_3s^3+\ua_4s^4+\\
 &\quad \oa_4s^5+\oa_6s^6+\ua_7s^7+\cdots, \\ %\nonumber
 K_3(s)&=\oa_0+\oa_1s+\ua_2s^2+\ua_3s^3+\oa_4s^4+\\
&\quad \oa_5s^5+\ua_6s^6+\ua_7s^7+\cdots, \\ %\nonumber
K_4(s)&=\oa_0+\ua_1s+\ua_2s^2+\oa_3s^3+\ua_4s^4+\\
&\quad \ua_4s^5+\ua_6s^6+\oa_7s^7+\cdots.  %\nonumber
\end{align*}
We have that 
if $K_1(s),\cdots,K_4(s)$ are Hurwitz then $P(s;[\ba])$ is Hurwitz. 
\begin{rem}
A trivial proof  of K4PT can be obtained by letting in (\ref{eq:K8}) $\ubb=0$ and $\obb=0$.
 We thus obtain $K_1(s)\equiv K_6(s)$, $K_2(s)\equiv K_5(s)$, $K_3(s)\equiv K_8(s)$, and $K_4(s)\equiv K_7(s)$
 and we are done.
 \end{rem}
\begin{rem}
\label{rem:conjH}
Here, since $\ba\in\RR^{N+1}$ we obtain: $P(s_0;\ba)=0$ iff $P(s_0^*;\ba)=0$, where  $s_0^*$ denotes the conjugate of $s_0$.
Hence, a containing sector for $P(s;[\ba])$ must have the symmetric form   $[0.5\pi+\alpha,1.5\pi-\alpha], \alpha\in(0,0.5\pi)$.
\end{rem}
In what follows we assume that $P(s;[\ba])$ is Hurwitz. 
By using Remark~\ref{rem:conjH}  we obtain that its minimal containing sector is $[0.5\pi+\alpha,1.5\pi-\alpha_0]$.
Hence,   we only need to compute  $\alpha_K\le\alpha_0$. Notice that   we can choose to rotate the roots of $P(s;[\ba])$ either CW or CCW
by $\alpha$.

We will rotate  the roots of $P(s;[\ba])$ CW by $\alpha$.  
The interval coefficients of the rotated polynomial  $P_\alpha(s;[\bc+i\bd])$ for $n=0,1,...,N$  are given by
\begin{align}
[c_n+id_n]&=e^{in\alpha}[a_n]\\                                                                   \nonumber
&=(\cos(n\alpha)+i\sin(n\alpha))[a_n]\\                                                          \nonumber
&=\cos(n\alpha)[a_n]+i\sin(n\alpha)[a_n],    \nonumber
\end{align}
where 
\begin{equation}
[c_n]=\cos(n\alpha)[a_n]
\end{equation}
and
\begin{equation}
[d_n]=\sin(n\alpha)[a_n]
\end{equation}
should be computed by using (\ref{eq:xy}).
Notice  the  dependence  of both $[c_n]$ and $[d_n]$ on $[a_n]$ and $\alpha$.

If there are no roots of  $P(s;[\ba])$ in the sector $(0.5\pi, 0.5\pi+\alpha]$
then    $P_\alpha(s;[\bc+i\bd])$ will be Hurwitz and we can check its stability by checking the stability of  its 8
Kharitonov polynomials. Then, by using the bisection algorithm we can compute $\alpha_K$. So, in order to decide
whether or not the sector
$[0.5\pi+\alpha,1.5\pi-\alpha]$ is a containing sector for $P(s;[\ba])$ we have to test the stability of a up to 8 complex Kharitonov polynomials
of  $P_\alpha(s;[\bc+i\bd])$.
\\ 
\\
\begin{example}
Let
\begin{align}
P(s;[\ba])=&[4.71,4.91]+[7.71,7.91]s+
       [3.9,4.1]s^2+[0.9,1.1]s^3. 
\end{align}
We derived this interval polynomial  by computing the coefficients of the polynomial whose roots are\\
 $s_1=-1.5+1.6i, s_2=-1.5-1.6i, s_3=-1$ and then
replaced each coefficient, say $a_n$, by\\
 $[a_n-0.1, a_n+0.1]$. 
The corresponding four Kharitonov polynomials are 
\begin{align}
K_1(s)&=4.71+7.71s+4.1s^2+1.1s^3\\ \nonumber
K_2(s)&=4.71+7.91s+4.1s^2+0.9s^3\\ \nonumber
K_3(s)&=4.91+7.91s+3.9s^2+0.9s^3 \\ \nonumber
K_4(s)&=4.91+7.71s+3.9s^2+1.1s^3.\\ \nonumber
\end{align}
These four Kharitonov polynomials are Hurwitz, therefore, $P(s;[\ba])$ is Hurwitz.
By combining  the proposed algorithm  and  the bisection algorithm we obtain that the containing sector ${\cal K}=[0.5\pi+\alpha_K,1.5\pi-\alpha_K]$ for  $P(s;[\ba])$ satisfies
$\alpha_K\in[\ualpha_K ,\oalpha_K]=[0.2005,0.2006]\pi=[36.09^\circ, 36.108^\circ]$.\\
Hence, ${\cal K}=[0.5\pi+\ualpha_K,1.5\pi-\ualpha_K]=[126.09^\circ,233.81^\circ]$ is a containing sector for $P(s;[\ba])$.\\
To check this result we generated a million uniformly distributed polynomials for $P(s;[\ba])$. These polynomials appeared in the sector
${\cal S}=[126.7957^\circ,  233.2043^\circ]$. We also checked all 16 vertex polynomials of $P(s;[\ba])$ and obtained that they are in the sector
${\cal V}=[126.7268^\circ,  233.2732^\circ]$ and were attained by
 the vertex polynomial $4.71+7.91s+3.9s^2+1.1s^3$.
Notice that ${\cal K}\supseteq {\cal O}\supseteq {\cal V}\supseteq {\cal S}$ and that we   were surprised  that ${\cal V}\supseteq {\cal S}$.
\end{example}
\begin{rem}
\label{rem:congR}
Based on the above example 
we conjecture that the minimal containing sector is attained at  a vertex polynomial of $P(s;[\ba])$.
If this is the case we have to test $2^{N+1}$  vertex polynomials unless we can prove that a smaller predetermined subset thereof suffices.
\end{rem}

\section{Conclusion}

Related algorithmic and interval-arithmetic techniques for robust stability analysis of N-dimensional discrete systems, including recursive branch-and-bound methods, are presented in \cite{HertzND}.

\label{sec:Conc}
In this article we  presented  algorithms to compute  a  suboptimal Kharitonov-type containing sectors, say 
${\cal K}=[0.5\pi+\alpha_K,1.5\pi-\beta_K]$, for
$P(s;[\ba+i\bb])$.\\
{\bf Complex case}: 
We gave an example for $N=3$ and computed  ${\cal K}$,
 the  sector ${\cal V}$ containing its $2^8$  vertex polynomials, and the sector ${\cal S}$  containing  a million polynomials
whose coefficients were chosen to be uniformly distributed in $[\ba+i\bb]$. We thus obtained that ${\cal K}\supseteq{\cal V} \supseteq {\cal S}$
and were surprised that ${\cal V}\supseteq{\cal S}$. 
Therefore, we conjectured that  ${\cal O}={\cal V}$, where ${\cal O}$ is the minimal containing sector. 
\\
{\bf Real case}: \\
We gave an example for $N=3$ and computed  ${\cal K}=[0.5\pi+\alpha_K,1.5\pi-\alpha_K]$,
the  sector ${\cal V}$ containing its $2^4$  vertex polynomials, and the sector ${\cal S}$  containing  a million polynomials
whose coefficients were chosen to be independently uniformly sampled  in each interval of $[\ba]$.We thus obtained that ${\cal K}\supseteq{\cal V} \supseteq {\cal S}$
and interestingly we obtained  that ${\cal V}\supseteq{\cal S}$. 
Therefore, we conjectured that  ${\cal O}={\cal V}$, where ${\cal O}$ is the minimal containing sector.

{\bf  Conjecture for complex case:} 
Suppose that the given $P(s;[\ba+i\bb])$  is Hurwitz then its minimal containing sector ${\cal O}[0.5\pi+\alpha_0,1.5\pi-\beta_0]$ is given by the  sector containing the   predetermined
set  $X$ of $P(s;[\ba+i\bb])$'s vertex polynomials, where $|X|\le 2^{2(N+1)}$.

{\bf  Conjecture for real case:} 
Suppose that the  given $P(s;[\ba])$  is Hurwitz then its minimal containing sector ${\cal O}[0.5\pi+\alpha_0,1.5\pi-\alpha_0]$ is given by the  sector containing the   predetermined
set  $X$ of $P(s;[\ba])$'s vertex polynomials, where $|X|\le 2^{N+1}$.\\

Identifying such a subset $X$ and establishing the conjecture formally remain interesting topics for future research.\\
 In addition, for the  Kharitonov-type containing sectors that has been computed for the real and complex examples, we carried out further tightening  by  the  recursive branch-and-bound interval-based robust stability methods \cite{HertzND}.  We could thus
 obtain  containing sectors that are very close to the vertrex containing sectors for both the real and complex cases.


\begin{thebibliography}{99}
\bibitem{RobCont} S.P. Bhattacharyya, H. Chapellat, and L.H. Keel,
\emph{ROBUST CONTROL The Parametric Approach,}  Prentice Hall PTR (1995).
\bibitem{KharitonovR}
V. L. Kharitonov, \emph{Asymptotic stability of an equilibrium position of a family of systems of differential equations}, Differentsialnye uravneniya, 14 (1978), 2086-2088. (in Russian).
\bibitem{Kharitonov8}
V. L. Kharitonov, \emph{On a Generalization of a Stability Criterion}, Akademii nauk Kazakhskoii  SSR,
Seria fisiko-mathematicheskaia,1 (1978): 53-57.
%\bibitem{KharitonovRQ} V.L. Kharitonov, \emph{The Routh-Hurwitz problem for families of polynomials and quasipolynomials???}, Izvetiy %Akademii Nauk Kazakhskoi SSR, Seria fizikomatematicheskaia, vol. 26 (1979), pp. 69-79.  %p627 in boo
\bibitem{HK} M.A. Hitz and E. Kaltofen, \emph{The Kharitonov   theorem and  its applications in symbolic mathematical computation}, Journal of Symbolic
Computation, 1997.
\bibitem{Bistriz} Y. Bistriz, \emph{Stability Criteria for Continuous-Time System Polynomials with Uncertain Complex Coefficients,}
IEEE Trans. on Circuits and Systems, Vol. 35, No. 4, Apr., 1988. pp. 442-448.
	\bibitem{HertzPreprint}
D.~Hertz, ``Algorithm to Compute a Kharitonov-Type Sector Containing All Roots
of a Real or Complex Hurwitz Interval Polynomial,''  Preprint, Sept.~2024.
doi:10.13140/RG.2.2.20452.39047, available online.
\bibitem{HertzND}
D.~Hertz,
``Robust Stability of Interval N-Dimensional Linear Discrete Systems
Possibly with Fractional Sampling via Recursive B\&B, Interval Arithmetic,
and Possibly Linear Optimization,''
Preprint, Jan.~2025.
doi:10.13140/RG.2.2.36041.66406, available online.

\end{thebibliography}
\end{document}